\documentclass[10pt]{amsart}
\usepackage[utf8]{inputenc}
\usepackage[margin=1in]{geometry}
\usepackage{hyperref}
\usepackage{amsmath}
\usepackage{amsfonts}
\usepackage{amssymb}
\usepackage{graphicx}
\usepackage{mathrsfs}
\usepackage{eqnarray,amsthm}
\usepackage{stmaryrd}
\usepackage{lmodern}
\usepackage[T1]{fontenc}
\usepackage{fullpage}
\usepackage{enumitem}
\usepackage{mathtools}
\usepackage{graphicx}
\usepackage{epsfig}
\usepackage{dsfont}
\usepackage{url}

\usepackage[usenames, dvipsnames]{color}

\newtheorem{thm}{\bf Theorem}[section]
\newtheorem{prop}[thm]{Proposition}
\newtheorem{defn}[thm]{Definition}

\newtheorem{lemma}[thm]{Lemma}
\newtheorem{corollary}[thm]{Corollary}
\newtheorem{theorem}[thm]{Theorem}
\theoremstyle{remark}
\newtheorem{remark}[thm]{Remark}

\newcommand\R{{\mathbb R}}
\newcommand{\cU}{\mathcal{U}}
\newcommand{\cDG}{\mathcal{DG}}
\newcommand{\dd}{{\, \mathrm d}}
\newcommand{\norm}[1]{\left\| #1 \right\|}

\begin{document}

\title{\textbf{Quantitative regularity for parabolic De Giorgi classes}}
%
%De Giorgi elliptic and parabolic regularity theory
\date{}

\author{Jessica Guerand}

%\address[Jessica Guerand]{University of Cambridge, Department of
%  Pure Mathematics and Mathematical Statistics, Wilberforce Road,
%  Cambridge CB3 0WA, United Kingdom}
%
%\email{jg900@cam.ac.uk}

\address[Jessica Guerand]{Universit\'e de Montpellier, IMAG, 499-554 Rue
  du Truel, 34090 Montpellier, France}
\email{jessica.guerand@umontpellier.fr}

\maketitle
\begin{abstract}
             We consider Parabolic De Giorgi classes for which we give a fully quantitative proof of the H\"older continuity for functions in these classes using the De Giorgi method. The main result is a quantitative proof of the intermediate value lemma, also known as the second lemma of De Giorgi. This paper is self-contained.  
\end{abstract}
          
\bigbreak
\bigbreak

\tableofcontents

          \section{Introduction}\label{intro}
          
This paper is an article version of the non-peer reviewed seminar paper \cite{gue3} about parabolic equations, here adapted to the more general framework of parabolic De Giorgi classes. Its goal is to give a short, light, and fully self-contained quantitative proof of the De Giorgi method for parabolic settings, which provides with alternative tools to existing results.

\subsection{The parabolic setting.}\label{Parabolic setting} Precisely here, the H\"older continuity result for the following parabolic equations in divergence form with rough coefficients, i.e. merely measurable coefficients and its corresponding De Giorgi classes, 
\begin{equation}
\label{eqpara}
\partial_{t}u = \nabla_{x}\cdot (A\nabla_{x} u) + B\cdot \nabla_{x} u+ S, \quad t\in \R \quad x\in\R^d,
\end{equation}
for $u=u(t,x)$, $A=A(t,x), B=B(t,x), S=S(t,x)$ satisfying 
\begin{equation}
\label{hypcoeff}
\begin{dcases}
\text{$A$ is a measurable $d\times d$ matrix field such that } 0<\lambda I\leq A \leq \Lambda I,\\
\text{$B$ is a measurable vector field of size $d$ such that } |B| \le \Lambda\\
\text{$S$ is a real scalar field in $L^{q}$ with } q>\max\left(2,\frac{d+2}{2}\right).
\end{dcases}
\end{equation}
 This class of equation \eqref{eqpara}-\eqref{hypcoeff} is invariant under translation $(t,x)\rightarrow (t+t_0,x+x_0)$ for any $z_0:=(t_0,x_0)\in \R\times \R^d$ and dilatation $(t,x)\rightarrow (r^2 t,r x)$ for $r>0$ of the variables but also any affine transformation of the solution $(t,x)\rightarrow \gamma u (t,x) + \rho$ is still a solution for any solution $u$ and any $(\gamma,\rho)\in \R^2$.  Following these invariances, define the parabolic cylinder $Q_{r}(z_0):= \left(t_0-r^2,t_0\right]\times B_{r}(x_0)$ where $B_r(x_0)$ is the euclidian ball of radius $r$ centered at $x_0$.  Introduce also the shorthand notation $Q_r:=\left(-r^2,0\right]\times B_{r}$ with $B_r$ the euclidian ball of radius $r$ centered at $0$.  We denote $Q_1^{-}=\left(-2,-1\right] \times B_1$.  We also define the positive (resp. negative) part of a function $u$ by $u_+=\max(u,0)$ (resp.  $u_-=\min(-u,0)$). The Lebesgue measure of a Lebesgue subset $E$ of $\mathbb{R}^N$ is denoted by $|E|$. 
 \vskip2mm
 
\begin{defn}[De Giorgi classes.]
\label{def:DGclass}
Let $\cU:=(T_1,T_2)\times \Omega$ be an open set of $\R\times \R^d$.  
 For non-negative parameters $\gamma_1,\gamma_2,\gamma_3$ and $p>0$, we define the De Giorgi sub-class (resp.  super-class) denoted by $\cDG^+(\gamma_1,\gamma_2,\gamma_3,p)$ (resp. $\cDG^-(\gamma_1,\gamma_2,\gamma_3,p)$) the set of function $u \in L^{\infty}((T_1,T_2);L^{2}(\Omega))\cap L^{2}((T_1,T_2),H^{1}(\Omega))$ which satisfies $\forall k\in \R$, $\forall (s,t) \in \R^2$ such that $T_1\le s<t\le T_2$, $\forall R>r>0$ and $\forall x_0 \in \Omega$ such that $B_R(x_0) \subset \Omega$, the following inequalities 
 \begin{align*}
 \int_{B_{r}(x_0)} (u-k)_{\pm}^{2} (t,x) \dd x + \gamma_1 \int_{s}^{t} \int_{B_{r}(x_0)} |\nabla_x (u-k)_{\pm}|^2(\tau,x) \dd x \dd \tau \\
\leq  \int_{B_{R}(x_0)} (u-k)_{\pm}^2 (s,x) \dd x + \frac{\gamma_2}{(R-r)^2} \int_{s}^{t} \int_{B_R(x_0)} (u-k)_{\pm}^{2} (\tau,x) \dd x \dd \tau \\ + \gamma_3 \left(\int_s^t \int_{B_R (x_0)} (u-k)_{\pm}^p(\tau,x)\dd x \dd \tau \right)^{1/p}.
 \end{align*}
\end{defn} 

\vskip2mm
\begin{remark}
The classes satisfy the same property of invariance than the equation: invariance of translation and dilation of the variables and any affine transformation of functions in classes (in particular multiplying by a negative constant send $\cDG^+$ to $\cDG^-$ and vice versa). 
\end{remark}

\subsection{Main contribution.} The invariance properties of the De Giorgi classes allow to state the theorems in unit cylinders. The two main contribution is a quantitative intermediate value lemma for Parabolic De Giorgi classes and the H\"older continuity theorem using a quantitative De Giorgi method.
 
 \vskip2mm
 \begin{theorem}[Intermediate value lemma for parabolic De Giorgi classes]
\label{t:LVI}
Let $\gamma_1, \gamma_2,\gamma_3, p >0$.  There is $C>0$ only depending on $d$, $\gamma_1, \gamma_2,\gamma_3, p$, such that any $u:Q_2\rightarrow\R \in \cDG^+(\gamma_1, \gamma_2,\gamma_3, p)$ such that $u\le 1$ in $Q_{\frac{3}{2}}$ satisfies 
\begin{equation}
\label{e:LVI}
|\{u\leq 0\}\cap Q_1^-||\{u\geq \frac{1}{2}\}\cap Q_1| \le C |\{0<u<\frac{1}{2}\}\cap Q_2|^{\frac{1}{4p+2}}.
\end{equation}
 \end{theorem}
 
 \vskip2mm
 \begin{remark}
 Cylinders $Q_1$ and $Q_1^-$ have to be disjoint and ordered in time this way. Indeed, the function $u(t,x)={\bf 1}_{t<-1}\in \cDG^+(\gamma_1, \gamma_2,\gamma_3, p)$ allows non-increasing jumps in time. 
 \end{remark}
  
  \vskip2mm
\begin{theorem}[H\"older continuity for parabolic De Giorgi classes]
\label{t:Holder}
Let $\gamma_1, \gamma_2,\gamma_3$ and $p\in \left[1,\min\left(2,\frac{d+2}{d}\right)\right)$.  
There is $\alpha \in (0,1)$ and $C>0$ only depending on $d$, $\gamma_1, \gamma_2,\gamma_3, p$, such that any $u : Q_2\rightarrow\R\in \cDG^+(\gamma_1, \gamma_2,\gamma_3, p)\cap \cDG^-(\gamma_1, \gamma_2,\gamma_3, p)$  satisfies 
\begin{equation*}
 [u]_{C^\alpha (Q_{1})} := \sup_{z_1,z_2 \in Q_1, \ z_1 \not =
      z_2} \frac{|u(z_1) - u(z_2)|}{|z_1-z_2|} \leq C \norm{u}_{L^2(Q_2)}.
\end{equation*}
\end{theorem} 

\subsection{Application to parabolic equations.}
Let us give the definition of weak solutions and the associated results.
 \vskip2mm
\begin{defn}[Weak sub/super-solutions]
Let $\cU:=(T_1,T_2)\times \Omega$ be an open set of $\R\times \R^d$.  
A function $u \hskip1mm:\hskip1mm \cU \rightarrow \R$ is a weak sub-solution (resp. super-solution) of \eqref{eqpara} if $$u\in L^{\infty}((T_1,T_2),L^{2}(\Omega))\cap L^{2}((T_1,T_2),H^1(\Omega)), \quad \partial_t u \in L^{2}((T_1,T_2),H^{-1}(\Omega))$$ and for all non-negative $\varphi \in~C^\infty_c (\cU)$, 
\begin{equation*}
-\int_{\cU} u\partial_t \varphi+ \int_{\cU} A\nabla_x u \cdot \nabla_x \varphi - \int_{\cU} B \cdot \nabla_x u \varphi -\int_{\cU} S\varphi \le 0 \quad \text{(resp. } \ge 0\text{)}.
\end{equation*}
A function $u \hskip1mm:\hskip1mm \cU \rightarrow \R$ is a weak solution of \eqref{eqpara} if $$u\in L^{\infty}((T_1,T_2),L^{2}(\Omega))\cap L^{2}((T_1,T_2),H^1(\Omega)),\quad \partial_t u \in L^{2}((T_1,T_2),H^{-1}(\Omega))$$ and for all $\varphi \in  C^\infty_c (\cU)$, 
\begin{equation*}
-\int_{\cU} u\partial_t \varphi+ \int_{\cU} A\nabla_x u \cdot \nabla_x \varphi - \int_{\cU} B \cdot \nabla_x u \varphi -\int_{\cU} S\varphi = 0.
\end{equation*}
\end{defn}
\vskip2mm
 \begin{corollary}[Intermediate value lemma for parabolic equation]
\label{cor:LVIsol}
There is $C>0$ only depending on $d$, $\lambda,  \Lambda$ and $q$, such that any weak sub-solution $u : Q_2\rightarrow\R$ of \eqref{eqpara}-\eqref{hypcoeff} such that $u\le 1$ in $Q_{\frac{3}{2}}$ satisfies 
\begin{equation}
\label{e:LVIsol}
|\{u\leq 0\}\cap Q_1^-||\{u\leq \frac{1}{2}\}\cap Q_1| \le C |\{0<u<\frac{1}{2}\}\cap Q_2|^{\frac{1}{4p+2}}.
\end{equation}
 \end{corollary}
 
 \vskip2mm
\begin{corollary}[H\"older continuity for parabolic equation]
\label{cor:Holder}
There is $\alpha \in (0,1)$ and $C>0$ only depending on $d$, $\lambda,  \Lambda$ and $q$ such that any weak solution $u:Q_2\rightarrow\R$ of \eqref{eqpara}-\eqref{hypcoeff} satisfies 
\begin{equation*}
 [f]_{C^\alpha (Q_{1})} := \sup_{z_1,z_2 \in Q_1, \ z_1 \not =
      z_2} \frac{|f(z_1) - f(z_2)|}{|z_1-z_2|} \leq C \norm{u}_{L^2(Q_2)}.
\end{equation*}
\end{corollary} 
 The proof of these corollary comes from the fact that sub/super-solutions are in De Giorgi classes. This proof can be found in Appendix \ref{A:sol in DG}.

\subsection{State of the art.} The De Giorgi Nash Moser techniques \cite{degiorgi,nash,moser} have been introduced to solve the $19^{\mathrm{th}}$ Hilbert problem about regularity of local minima of functionals. The gradients of these minima are weak-solutions to Euler-Lagrange equations. Seeing these quasilinear elliptic equations as linear ones in divergence form with merely measurable coefficients allows to get the H\"older continuity of solutions. 
Concerning solutions of parabolic equations, we refer to \cite{gue3} for an historical overview.

Initially formulated by Ladyzhenskaya and Uralt'seva \cite{LU1} in the elliptic case and in \cite{LSU} in the parabolic case, De Giorgi classes were used to deal not only with solutions of equations with rough coefficient but also local minima and quasi-minima of functional that do not necessarily satisfy any Euler-Lagrange equations, see Giaquinta and Gusti \cite{Giaquinta_Gusti82} for the elliptic case and Wieser \cite{Wieser87} for the parabolic case. Several works followed. Let us mention the work of DiBenedetto \cite{Dibenedetto89} with a Harnack inequality for De Giorgi classes without the use of any covering argument. Also the work of Gianazza and Vespri \cite{Gianazza_Vespri06} deals with De Giorgi classes of order $p>1$. Later on Masson and Silijander \cite{Masson_Silijander13} extended the results to metric spaces for $p\geq 2$. Meanwhile the first version of this paper, Liao \cite{Liao21} wrote remarks on Harnack inequalities for parabolic De Giorgi classes with new perspectives and later Her\'an \cite{Heran22} proved H\"older continuity of parabolic quasiminimiser in metric spaces. All theses techniques are in the spirit of Moser's ideas. 

In this paper, the Giorgi techniques are used to give a new short quantitative self-contained and alternative way to deal with parabolic De Giorgi classes.  
The main contribution is the quantitative intermediate value lemma in a parabolic framework, also called second lemma of De Giorgi or isoperimetric inequality. This lemma is quantitative in elliptic framework \cite{degiorgi,Vasseur2,HouNiu,DiBen2010} using either isoperimetric arguments or a Poincar\'e inequality.
About parabolic framework, no quantitative version of this lemma seems to exist in the local case. One can find non-quantitative versions, for example in \cite{Vasseur2}, a version obtained by contradiction with a compactness argument. There exists a quantitative version of this lemma for nonlocal time-dependent integral operator \cite{CCV} which does not seem to apply for either local parabolic equations nor parabolic De Giorgi classes.
In this paper, we introduce a new point of view, using measure properties directly from the definition of the classes. We also provide an alternative way for proving the $L^{\infty}$ bound without using the interpolation inequality in \cite{Vasseur2}.

\subsection{Strategy of the proof.}
The core of the proof, which is the intermediate value lemma, also called the second lemma of De Giorgi, is introduced in Section \ref{Sec:IVL}. The $L^{\infty}$ bound, sometimes called first lemma of De Giorgi, is given in Section \ref{Sec:DGmethod}. Then it is shortly recalled how to get the Measure to Pointwise lemma from both De Giorgi's lemma and from it how to get the H\"older continuity thanks to the oscillation.

\section{Intermediate value lemma}
\label{Sec:IVL}
 
 Before giving the result for functions in $\cDG^+$, we recall it for functions in $H^1$.
 \subsection{Functions in $H^1$.} Let us state the following result that we use to get the intermediate value lemma for functions in $\cDG^+$. We give a different version of the so-called isoperimetric inequality \cite[Lemma 10]{Vasseur2} and an alternative proof of the paper \cite[Theorem 2.9]{HouNiu}.
 
 \begin{lemma}[Intermediate value lemma in $H^1$]
 \label{l:IVL H1}
 There is $C(d)>0$ only depending on the dimension such that any $u \in H^1(B_1)$ satisfies
 \begin{equation*}
 |\{u\le 0\}\cap B_1||\{u\ge \frac{1}{2}\}\cap B_1| \le C(d) |\{0<u < \frac{1}{2}\}\cap B_1|^{\frac{1}{2}} \left(\int_{B_1} |\nabla u_+|^2 \right)^{\frac{1}{2}}.
 \end{equation*}
 \end{lemma}
 
 \begin{proof}
 Apply a Poincar\'e inequality to the truncated function 
 \begin{equation}
v(x):= \left\{ \begin{array}{lll}
 0 & \mbox{ if } & u(x)\le 0,\\ 
 u(x) & \mbox{ if } & 0<u(x)<\frac{1}{2}, \\
 \frac{1}{2} & \mbox{ if } & u(x) \ge \frac{1}{2}.
 \end{array}\right. \in H^{1}(B_1) \subset W^{1,1}(B_1), 
 \end{equation}
 see for example \cite[Theorem 3.2]{AD04},
 \begin{align}
 \int_{B_1} \left|v(x)-\int_{B_1} v(y) \dd y\right| \dd x \le C(d) \int_{B_1} |\nabla v(x)| \dd x.
 \end{align}

Computing a bound from below for the left hand side gives 
 \begin{align}
\label{ineq1}
 \displaystyle \int_{B_{1}} \left|v(x)-\int_{B_1} v(y) \dd y\right| \, \mathrm{d}x \geq \displaystyle  \int_{\{v=0 \}\cap B_1} \left|v(x)-\int_{B_1} v(y) \dd y\right| \, \mathrm{d}x \geq \displaystyle \int_{\{v=0 \}\cap B_1}\int_{B_1} v(y) \dd y \, \mathrm{d}x  \nonumber \\ \geq \frac{1}{2|B_{1}|}|\{v=0\}\cap B_1| |\{v=\frac{1}{2}\}\cap B_1 | =\frac{1}{2|B_{1}|}|\{u\leq 0\}\cap B_1| |\{u\geq \frac{1}{2}\}\cap B_1 |
 \end{align}
 
 and from above for the right hand side using a Cauchy-Schwarz inequality gives
 \begin{align}
 \label{ineq2}
 \displaystyle \int_{B_{1}} | \nabla v (x) | \, \mathrm{d}x = \int_{\{0<u<\frac{1}{2}\}\cap B_1} | \nabla v (x) | \, \mathrm{d}x \leq \sqrt{\int_{B_{1}} | \nabla (u-k)_{+} (x) |^{2}  \, \mathrm{d}x}  |\{k<u<l\}\cap B_1|^{\frac{1}{2}}.
 \end{align}
 Combining \eqref{ineq1} and \eqref{ineq2} gives the result.
 \end{proof}
 
 \subsection{Functions in $\cDG^+$.} To prove the intermediate value lemma , Theorem \ref{t:LVI}, which is the main novel proof of the paper, two lemmas are required. Notice that the left hand side of the De Giorgi class inequalities has two non-negative terms so that it can be split in two inequalities. This first following lemma comes from the second part of the inequality in Definition~\ref{def:DGclass} which is the Caccioppoli like inequality.
 
  \vskip2mm
\begin{lemma}[Universal bound of the $L^2$ of the gradient]
\label{Lm:LVI elliptic}
Let $u:Q_{2}\rightarrow \mathbb{R}$ be a function in $DG^{+}(\gamma_1,\gamma_2,\gamma_3,p)$ such that $u\leq 1$ on $Q_{\frac{3}{2}}$.
Then there exists a constant $C>0$ which only depends on $d, \gamma_1, \gamma_2, \gamma_3$ and $p$ such that for all $-2<s<t<0$, we have
 $$\displaystyle \int_{s}^{t}\int_{B_{\frac{5}{4}}} |\nabla_{x}u_+|^{2}(\tau,x) \, \mathrm{d}x\mathrm{d}\tau \leq C.$$
\end{lemma}

\begin{proof}[Proof of Lemma \ref{Lm:LVI elliptic}]
Use Definition \ref{def:DGclass} for $r=\frac{5}{4}$, $R=\frac{3}{2}$, $x_0=0$ and deduce 
$$ \displaystyle \int_{s}^{t}\int_{B_{\frac{5}{4}}} |\nabla_{x}u_+|^{2}(\tau,x) \, \mathrm{d}x\mathrm{d}\tau \leq \frac{|B_{\frac{3}{2}}|}{\gamma_1}+32\frac{\gamma_2}{\gamma_1}|B_{\frac{3}{2}}|+2^{\frac{1}{p}}\frac{\gamma_3}{\gamma_1}|B_{\frac{3}{2}}|^{\frac{1}{p}},$$
which gives the result.
\end{proof} 
 
This second following lemma comes from the first part of the inequality in Definition~\ref{def:DGclass}.

\begin{lemma}[A key inequality for close times]
\label{Lm2LVIP}
Let $u:Q_{2}\rightarrow \mathbb{R}$ be a function in $DG^{+}(\gamma_1,\gamma_2,\gamma_3,p)$ such that $u\leq 1$ on $Q_{\frac{3}{2}}$.Then for all $(t_{1},t_{2}, \tau)\in (-2,0)^3$ such that $-2<t_1<\tau<t_2<0$, we have
   \begin{equation*}
  |u\geq 1/2,(\tau,t_2)\times B_1||u\leq 0, (t_1,\tau)\times B_1| \leq C |0<u<1/2, (t_1,\tau)\times B_{2}|^{\frac{1}{2}}\nonumber\\ 
 + C (t_2-t_1)^{2+\frac{1}{p}},
 \end{equation*}
 with $C$ only depending on $d, \gamma_1, \gamma_2, \gamma_3$ and $p$.
\end{lemma}

\vskip2mm
\begin{proof}[Proof of Lemma \ref{Lm2LVIP}]
In this proof, let $C>0$ be a constant which only depends on $d, \gamma_1, \gamma_2, \gamma_3$ and $p$ which may change at each line.
Thanks to Definition \ref{def:DGclass}, using that $u\leq 1$ and that $(t-s)\leq C(t-s)^{1/p}$
\begin{align}
\label{ineq no jump time}
 \displaystyle \int_{B_1} u_{+}^{2}(t,x)\mathrm{d}x 
  \leq  \int_{B_{\frac{5}{4}}} u_{+}^2
 (s,x) \mathrm{d}x
 + C(t-s)^{1/p}.
 \end{align}
Compute a bound from below of the left hand side 
 \begin{align}
 \label{bd lhs}
 \displaystyle \int_{B_1} u_{+}^{2}(t,x)\mathrm{d}x\geq  \displaystyle \int_{\{y\in B_1, u(t,y)\geq 1/2 \}} (1/2)^{2}\mathrm{d}x \geq C |\{u(t,.)\geq 1/2\}\cap B_1|.
 \end{align}
 Compute a bound from below of the first term of the right hand side 
 \begin{align}
 \label{bd rhs 1}
 \displaystyle \int_{B_{\frac{5}{4}}} u_{+}^2
 (s,x) \mathrm{d}x & \leq \int_{\{y\in B_{\frac{5}{4}}, 0 < u(s,y) < 1/2 \}} u_{+}^2
 (s,x) \mathrm{d}x +  \int_{\{y\in B_{\frac{5}{4}}, u(s,y)\geq 1/2 \}} u_{+}^2
 (s,x) \mathrm{d}x \nonumber\\
 & \leq C\left(|\{0<u(s,.)<1/2 \}\cap B_{\frac{5}{4}}| + |\{u(s,.)\geq 1/2 \}\cap B_{\frac{5}{4}}| \right).
 \end{align}
 So combining \eqref{ineq no jump time}, \eqref{bd lhs} and \eqref{bd rhs 1} gives
 \begin{align*}
 |\{u(t,.)\geq 1/2\}\cap B_1| \leq &  C \left(|\{0<u(s,.)<1/2 \}\cap B_{\frac{5}{4}}| + |\{u(s,.)\geq 1/2 \}\cap B_{\frac{5}{4}}| \right) + C(t-s)^{\frac{1}{p}}.
 \end{align*}
 
Multiply the last inequality by $|\{u(s,.)\leq 0 \}\cap B_{\frac{5}{4}}|$ and use Lemma \ref{l:IVL H1} and \ref{Lm:LVI elliptic} for $u(s,.)\in H^{1}(B_{\frac{5}{4}})$ by Fubini's theorem, 
\begin{align*}
%\label{concl1}
& |\{u(t,.)\geq 1/2\}\cap B_1| |\{u(s,.)\leq 0 \}\cap B_{1}|\\
  \leq & C \bigg(|\{0<u(s,.)<1/2 \}\cap B_{\frac{5}{4}}| 
  +|\{0<u(s,.)<1/2 \}\cap B_{\frac{5}{4}}|^{\frac{1}{2}}\sqrt{\displaystyle \int_{B_{\frac{5}{4}}} |\nabla_{x}u_+|^{2}(s,x) \, \mathrm{d}x} \bigg) + C (t-s)^{\frac{1}{p}}.
 \end{align*}
Integrate the latter over $s\in (t_{1},\tau)$ and $t\in (\tau,t_{2})$ with $-2\leq t_{1}<\tau<t_{2}\leq 0$ and obtain using a Cauchy-Schwarz inequality and Lemma \ref{Lm:LVI elliptic},
 \begin{align}
 \label{concl2}
   &|\{u\geq 1/2\}\cap (\tau,t_2)\times B_1||\{u\leq 0\}\cap (t_1,\tau)\times B_1| \nonumber\\
 \leq & C \bigg(|\{0<u<1/2\} \cap (t_1,\tau)\times B_{\frac{5}{4}}|
 +|\{0<u<1/2\}\cap (t_1,\tau)\times B_{\frac{5}{4}}|^{\frac{1}{2}}\bigg) + C (t_2-t_1)^{2+\frac{1}{p}}\nonumber\\
 \leq &  C |\{0<u<1/2\}\cap (t_1,\tau)\times B_{2}|^{\frac{1}{2}} + C (t_2-t_1)^{2+\frac{1}{p}}.
 \end{align}
 which ends the proof.
\end{proof}

Now we give a proof of Theorem \ref{t:LVI}. The idea is to understand that the ``error'' term $(t_2-t_1)^{2+\frac{1}{p}}$ in Lemma \ref{Lm2LVIP} is negligible compared to the other terms when $t_2-t_1$ is small and when the intervals are well-chosen.

\begin{proof}[Proof of Theorem \ref{t:LVI}]
Noticing that $|\{u\leq 0\}\cap (t_1,\tau)\times B_1|= |\{u< 1/2\}\cap (t_1,\tau)\times B_1|- |\{0<u<1/2\}\cap (t_1,\tau)\times B_1|$, Lemma \ref{Lm2LVIP} gives also the following inequality
\begin{align}
\label{ineq avec l}
 |\{u\geq 1/2\}\cap (\tau,t_2)\times B_1||\{u< 1/2\}\cap (t_1,\tau)\times B_1| \leq  C |\{0<u<1/2\}\cap Q_2|^{\frac{1}{2}}
  + C (t_2-t_1)^{2+\frac{1}{p}}.
\end{align}
We discretize the time interval in intervals of size $\frac{1}{n}$ for $n\in \mathbb{N}\setminus \{0\}$ and define for $k\in \llbracket 0,2n\rrbracket$,  $t_{k}=-2+\frac{k}{n}$. For $m\in \llbracket 1,2n-1\rrbracket$ let $I_m$ be the inequality 
 \begin{equation}
 \label{const1}
 |\{u< 1/2\}\cap (t_{m-1},t_{m})\times B_1|  \geq \frac{|u< 1/2 \mbox{, } Q_1^-|}{2n},
 \end{equation}
 and $J_m$ the inequality 
  \begin{equation}
 \label{const2}
 |\{u\geq 1/2\}\cap (t_{m},t_{m+1})\times B_1| \geq \frac{|u\geq 1/2 \mbox{, } Q_1|}{2n}.
 \end{equation}
By the pigeonhole principle, there exists $i\in \llbracket 1,n\rrbracket$ such that $I_i$ holds and $j\in \llbracket n, 2n-1\rrbracket$ such that $J_j$ holds. 
Notice that if $I_{m+1}$ doesn't hold then $J_m$ holds. Indeed if
\begin{align*}
|\{u< 1/2\}\cap (t_{m},t_{m+1})\times B_1|< \frac{|\{u< 1/2\} \cap Q_1|}{2n},
\end{align*}
then
\begin{align*}
 |\{u\geq 1/2\} \cap (t_{m},t_{m+1})\times B_1| \geq \frac{|B_{1}|}{n}-\frac{|\{u< 1/2\} \cap Q_1|}{2n}\geq \frac{|\{u\geq 1/2\}\cap Q_1^-|}{2n}.
\end{align*}
So necessarily there is $p\in \llbracket i,j\rrbracket$ such that $I_p$ and $J_p$ hold. Then use \eqref{ineq avec l} to get
 $$\begin{array}{lll}
\frac{|\{u<1/2\}\cap Q_1^-|}{2n} \frac{|\{u\geq 1/2\} \cap Q_1|}{2n} & \leq & |\{u< 1/2\}\cap (t_{p-1},t_{p})\times B_1| |\{u \geq  1/2\} \cap (t_{p},t_{p+1})\times B_1|\\
 & \leq & |\{0<u<1/2\}\cap Q_2|^{\frac{1}{2}}+C\left(\frac{2}{n}\right)^{2+\frac{1}{p}},
 \end{array}$$
  which leads to
  $$ |\{u< 1/2\} \cap Q_1^-| |\{u\geq 1/2\} \cap Q_1|\leq C n^2|\{0<u<1/2\}\cap Q_2|^{\frac{1}{2}}+C n^{-\frac{1}{p}}.$$
 So necessarily $|\{0<u<1/2\}\cap Q_2|>0$. And taking $n$ such that $$Cn^{-\frac{1}{p}}\leq C \frac{n^2|\{0<u<1/2\}\cap Q_2|^{\frac{1}{2}}}{2},$$ for example $n=\Big\lfloor \frac{2}{|\{0<u<1/2\}\cap Q_2|^{\frac{p}{4p+2}}} \Big\rfloor +1$, we get 
 $$|\{u< 1/2\}\cap Q_1^-| |\{u\geq 1/2\}\cap Q_1|\leq C|\{0<u<1/2\}\cap Q_2|^{\frac{1}{4p+2}}$$
 which achieves the proof of the Theorem since $|\{u< 1/2\}\cap Q_1^-| \geq |\{u\leq 0\}\cap Q_1^-|$.
\end{proof}
\section{De Giorgi method for parabolic De Giorgi classes}
\label{Sec:DGmethod}
This section is a short recall of the other steps of the De Giorgi method from \cite{Vasseur2} adapted to the parabolic De Giorgi classes. Only the next lemma uses a slightly different technique since it does not involve an interpolation inequality anymore \cite{Vasseur2}.  

\subsection{First Lemma of De Giorgi.} Let us state the first lemma of De Giorgi. 

\vskip2mm
\begin{lemma}[First Lemma of De Giorgi: $L^2-L^{\infty}$ estimate]
\label{firstlm}
Let $\gamma_1, \gamma_2,\gamma_3$ and $p\in \left[1,\min\left(2,\frac{d+2}{d}\right)\right)$.  
There exists a positive constant $\delta$ which only depends on $d, \gamma_1$, $\gamma_2, \gamma_3$ and $p$  such that 
for any $u:Q_{2}\rightarrow \mathbb{R}$ in $DG^{+}(\gamma_1,\gamma_2,\gamma_3,p)$ the following implication holds true.
If 
$$\displaystyle \int_{Q_{1}} u_{+}^2 \leq \delta,$$
then we have 
$$u_{+}\leq \frac{1}{2} \quad \mbox{ in } Q_{1/2}.$$
\end{lemma}

This proof is an alternative version of the ones in \cite{Vasseur2, GIMV} which does not involve an interpolation inequality and no $L^p$ gain of integrability but still involves the sequential inequality. 

\begin{proof}[Proof of Theorem \ref{firstlm}]
In this proof $C>0$ will denote a constant which will only depend on $d, \gamma_1, \gamma_2, \gamma_3,$ and $p$. 
Define 
$$U_k=\displaystyle\int_{Q_{r_k}} (u-c_k)_+^2 \mathrm{d}x\mathrm{d}t, $$
with $r_k = \frac{1}{2}(1+2^{-k})  $ and $c_k = \frac{1}{2}(1-2^{-k}) $. 
%We notice that $Q_{r_k}$ goes from $Q_1$ to $Q_{\frac{1}{2}}$ and $c_k$ from $0$ to $\frac{1}{2}$.    
Defining $V_k=U_{2k}$, the idea is to prove that the sequence $(V_k)$ satisfies
\begin{equation}
\label{e:induction}
V_k\leq C^{k} V_{k-1}^{\lambda},
\end{equation}
and deduce that $V_k=U_{2k}$ tends to $0$ when $V_0=U_0$ is small enough thanks to Lemma \ref{seq tends 0} in Appendix \ref{B:tool}. Since $U_0 = \displaystyle \int_{Q_{1}} u_{+}^2$ and $U_{\infty}=  \displaystyle \int_{Q_{\frac{1}{2}}} \left(u-\frac{1}{2}\right)_{+}^2=0$, the result follows. 
Let us prove the induction formula \eqref{e:induction}. 
Define the exponent  $$\rho= \left\{ \begin{array}{ll}
\frac{2d}{d-2} &\mbox{ if } d>2, \\
\frac{3}{2-p} & \mbox{ if } d \in \{1,2\}
\end{array}\right.$$ 
 in the following Sobolev inequality, for almost every $t\in (-r_k^2,0),$
\begin{equation}
\label{sobolevineq}
\| (u-c_k)_+(t,\cdot) \|_{L^\rho (B_{r_{k}})} \leq C(d)\| (u-c_k)_+(t,\cdot) \|_{H^{1}(B_{r_{k}})},
\end{equation}
where $C(d)$ is a constant which only depends on the dimension $d$ explicit in \cite{adams03}.
Using an H\"older inequality, get
\begin{align}
\label{lm1 ineq1}
U_k=\displaystyle\int_{Q_{r_k}} (u-c_k)_+^2 
\leq \displaystyle\int_{-r_k^2}^{0} \left(\displaystyle\int_{B_{r_{k}}} (u-c_k)_+^\rho (t,\cdot) \mathrm{d}x\right)^{\frac{2}{\rho}} |\{u(t,\cdot)\geq c_k\}\cap B_{r_k}|^{1-\frac{2}{\rho}} \mathrm{d}t.
\end{align}
Since $\{u(t,\cdot)\geq c_k \} = \{u(t,\cdot) \geq c_{k-1}+ 2^{-k-1} \}$, deduce that 
\begin{align}
\label{lm1 ineq2}
 |\{u(t,\cdot)\geq c_k\}\cap B_{r_k}|^{1-\frac{2}{\rho}} &\leq  |\{u(t,\cdot)\geq c_{k-1}+ 2^{-k-1} \}\cap B_{r_k}|^{1-\frac{2}{\rho}} \nonumber \\
 & \leq \left(2^{2k+2} \displaystyle\int_{B_{r_{k}}} (u-c_{k-1})_+^2(t,\cdot) \right)^{1-\frac{2}{\rho}} \nonumber \\
 &\leq C^k \left( \sup\limits_{t\in (-r_k^2,0)} \displaystyle\int_{B_{r_{k}}} (u-c_{k-1})_+^2(t,\cdot) \right)^{1-\frac{2}{\rho}}.
\end{align}
Use the first part of the inequality in Definition \ref{def:DGclass} with $s$ integrated in $(-r_{k-1}^2,-r_k^2)$ to bound the supremum in \eqref{lm1 ineq2} and $(u-c_{k-1})_+^p\leq  (u-c_{k-1})_+^2+ \mathds{1}_{\{u\geq c_{k-1} \}}$,
 \begin{align}
\label{lm1 ineq3}
 |\{u(t,\cdot)\geq c_k\}\cap B_{r_k}|^{1-\frac{2}{\rho}} 
 &\leq  C^k \left(\displaystyle\int_{-r_{k-1}^2}^{0} \int_{B_{r_{k-1}}} (u-c_{k-1})_{+}^{2} +  \left(\displaystyle\int_{-r_{k-1}^2}^{0} \int_{B_{r_{k-1}}} (u-c_{k-1})_{+}^p\right)^{\frac{1}{p}} \right)^{1-\frac{2}{\rho}} \nonumber \\
 &\leq   C^k \left(\displaystyle\int_{-r_{k-1}^2}^{0} \int_{B_{r_{k-1}}} (u-c_{k-1})_{+}^{2} +  \displaystyle\int_{-r_{k-1}^2}^{0} \int_{B_{r_{k-1}}} \mathds{1}_{\{u\geq c_{k-1} \}} \right)^{1-\frac{2}{\rho}} \nonumber\\
 &+ C^k \left(\displaystyle\int_{-r_{k-1}^2}^{0} \int_{B_{r_{k-1}}} (u-c_{k-1})_{+}^{2} +  \displaystyle\int_{-r_{k-1}^2}^{0} \int_{B_{r_{k-1}}} \mathds{1}_{\{u\geq c_{k-1} \}} \right)^{\frac{1}{p}(1-\frac{2}{\rho})}.
\end{align}
 %&  \leq  C^k \left(U_{k-1} + U_{k-1}^{\frac{1}{2}} \right)^{1-\frac{2}{p}},
Since 
\begin{align*}
  \displaystyle\int_{-r_{k-1}^2}^{0} \int_{B_{r_{k-1}}} \mathds{1}_{\{u\geq c_{k-1} \}} & = |\{u\geq c_{k-1} \}\cap Q_{r_{k-1}} |\\
 & \leq |\{u(t,\cdot)\geq c_{k-2}+ 2^{-k} \}\cap Q_{r_k-1}| \\
& \leq 2^{2k} \displaystyle\int_{Q_{r_{k-1}}} (u-c_{k-2})_+^2 \\
& \leq 2^{2k} U_{k-2},
\end{align*}
deduce using \eqref{lm1 ineq3}, 
\begin{align}
\label{lm1 ineq4}
 |\{u(t,\cdot)\geq c_k\}\cap B_{r_k}|^{1-\frac{2}{\rho}} & \leq C^{k} \left( \left( U_{k-1} + U_{k-2} \right)^{1-\frac{2}{\rho}}+\left( U_{k-1} + U_{k-2} \right)^{\frac{1}{p}(1-\frac{2}{\rho})} \right)\nonumber\\
 & \leq  C^{k} \left(U_{k-2}^{1-\frac{2}{\rho}} + U_{k-2}^{\frac{1}{p}(1-\frac{2}{\rho})}\right)
\end{align}
Notice that the last bound is independent of the variable $t$ so it remains to bound $\displaystyle\int_{-r_k^2}^{0} \left(\displaystyle\int_{B_{r_{k}}} (u-c_k)_+^\rho (t,\cdot) \mathrm{d}x\right)^{\frac{2}{\rho}} \mathrm{d}t$ in \eqref{lm1 ineq1}.
Using the Sobolev inequality \eqref{sobolevineq} and the second part of Definition \ref{def:DGclass} with $s$ integrated in $(-r_{k-1}^2,-r_k^2)$, deduce
 \begin{align}
 \label{lm1 ineq5}
\displaystyle\int_{-r_k^2}^{0} \left(\displaystyle\int_{B_{r_{k}}} (u-c_k)_+^\rho (t,\cdot)\right)^{\frac{2}{\rho}}
& \leq   C\left(\displaystyle\int_{-r_k^2}^{0}\displaystyle\int_{B_{r_{k}}} (u-c_k)_+^2(t,\cdot) +\displaystyle\int_{-r_k^2}^{0} \displaystyle\int_{B_{r_{k}}} |\nabla_x (u-c_k)_+|^2(t,\cdot) \right) \nonumber \\
 &\leq   C \left(\displaystyle\int_{-r_{k-1}^2}^{0} \int_{B_{r_{k-1}}} (u-c_{k-1})_{+}^{2} +  \left(\displaystyle\int_{-r_{k-1}^2}^{0} \int_{B_{r_{k-1}}} (u-c_{k-1})_{+}^p\right)^{\frac{1}{p}} \right)\nonumber \\
&\leq  C \left(\displaystyle\int_{-r_{k-1}^2}^{0} \int_{B_{r_{k-1}}} (u-c_{k-1})_{+}^{2} +  \displaystyle\int_{-r_{k-1}^2}^{0} \int_{B_{r_{k-1}}} \mathds{1}_{\{u\geq c_{k-1} \}} \right)\nonumber\\
 &+ C  \left(\displaystyle\int_{-r_{k-1}^2}^{0} \int_{B_{r_{k-1}}} (u-c_{k-1})_{+}^{2} +  \displaystyle\int_{-r_{k-1}^2}^{0} \int_{B_{r_{k-1}}} \mathds{1}_{\{u\geq c_{k-1} \}} \right)^{\frac{1}{p}}\nonumber\\
 &\leq   C^k \left(U_{k-2}+ U_{k-2}^{\frac{1}{p}}\right).
 \end{align}
By definition $U_k$ is non-increasing so assuming that $U_0<1$ gives $U_k<1$ for every $k\geq 0$. 
Combining \eqref{lm1 ineq4} and \eqref{lm1 ineq5} in \eqref{lm1 ineq1} and assuming $U_0<1$, deduce that $U_k$ satisfies the formula
$$U_k\leq C^k \left(U_{k-2}+ U_{k-2}^{\frac{1}{p}}\right) \left(U_{k-2}^{1-\frac{2}{\rho}} + U_{k-2}^{\frac{1}{p}(1-\frac{2}{\rho})}\right) \leq C^{k} U_{k-2}^{\alpha},$$
with $\alpha=\frac{1}{p}\left(2-\frac{2}{\rho}\right)>1$ which ends the proof using Lemma \ref{seq tends 0} choosing $\delta<C^{-\frac{\alpha^2}{(\alpha-1)^2}}$.
\end{proof}

\subsection{Measure to Pointwise Lemma.} The following Lemma requires the use of a sequence of functions that should remains in the De Giorgi class $DG^{+}(\gamma_1,\gamma_2,\beta\gamma_3,p)$ up to a specific index. This is satisfied starting with a function $v$ in $DG^{+}(\gamma_1,\gamma_2,\beta\gamma_3,p)$, with $\beta$ given by 
\begin{equation}
\label{beta}
\beta:=\frac{1}{2^{|Q_2|\left(\frac{2C}{\delta|Q_1|}\right)^{4p+2}+1}}.
\end{equation}
 with $\delta$ given in Lemma \ref{firstlm} and $C$ given in Lemma \ref{t:LVI}.

 \vskip2mm
\begin{lemma}[Measure to Pointwise]
\label{DGpLm5}
There exists a constant $\mu\in (0,1)$ which only depends on $d, \gamma_1, \gamma_2, \gamma_3,$ and $p$, such that for any function $v:Q_{2}\rightarrow \mathbb{R}$  in $DG^{+}(\gamma_1,\gamma_2,\beta\gamma_3,p)$ where $\beta$ is given in \eqref{beta}, if $v$ verifies
\begin{equation}
\label{hypvplm5}
\left\{ \begin{array}{c}
v\leq 1 \mbox{ in } Q_{\frac{3}{2}}\\
|\{v\leq 0 \}\cap \overline{Q}_1| \geq \frac{|\overline{Q}_1|}{2},
\end{array}\right.
\end{equation}
then
$$v\leq 1-\mu \quad \mbox{ in } Q_{\frac{1}{2}}.$$
\end{lemma}

\begin{proof}[Proof of Lemma \ref{DGpLm5}]
Introduce a sequence of function $v_k$ remaining in $DG^{+}(\gamma_1,\gamma_2,\gamma_3,p)$ for $k\leq k_0:=|Q_2|\left(\frac{2C}{\delta|Q_1|}\right)^{4p+2}+1$,
$$\left\{\begin{array}{l}
v_{0}=v\\
v_{k}=2\left(v_{k-1}-\frac{1}{2}\right). 
\end{array} 
\right.$$ 
The sequence satifies $v_{k}=2^k\left(v-(1-2^{-k}) \right)$, the sets $\{ 0<v_{k}<\frac{1}{2}\}= \{1-\frac{1}{2^{k}}<v<1-\frac{1}{2^{k+1}}\}$ are disjoint and the sequence $v_k$ still satisfies \eqref{hypvplm5}.

Necessarily, there exists $k_1\leq k_0$ such that $\displaystyle \int_{Q_{1}} (v_{k_1})_{+}^2 \leq \delta$.
Indeed, the intermediate value lemma, Lemma \ref{t:LVI} gives for any $k\leq k_0$ 
\begin{align*}
C|\{ 0<v_{k}< \frac{1}{2}\}\cap Q_{2}|^{\frac{1}{4p+2}} &\geq |\{v_{k}\leq 0\} \cap Q_1^-| |\{v_{k}\geq \frac{1}{2}\} \cap Q_{1}|\\
&\geq |\{v\leq 0 \}\cap Q_1^-||\{v_{k+1}\geq 0\} \cap Q_{1}|\\
&\geq \frac{|Q_1|}{2} \displaystyle \int_{Q_{1}} (v_{k+1})_{+}^2.
\end{align*}
Summing all the intermediate measure from $k=n+1$ to $n+m$ for any $n\geq 0$ and $m\geq 1$ such that $\int_{Q_{1}} (v_{k+1})_{+}^2 >\delta$ for $n+1\leq k \leq n+m$ and using the fact that the sets are disjoint gives 
\begin{align*}
|Q_2|\geq \sum_{k=n+1}^{n+m} | \{ 0<v_{k}< \frac{1}{2}\}\cap Q_{2}|\geq  \left(\frac{|Q_1|}{2C} \right)^{4p+2}\sum_{k=n+1}^{n+m} \left( \displaystyle \int_{Q_{1}} (v_{k+1})_{+}^2\right)^{4p+2} 
\geq  m\left(\frac{\delta|Q_1|}{2C} \right)^{4p+2},
\end{align*}
so that $m\leq |Q_2|\left(\frac{2C}{\delta|Q_1|}\right)^{4p+2}<k_0$ which means there is no interval of $\mathbb{N}$ of size $k_0$ such that $\int_{Q_{1}} (v_{k+1})_{+}^2 >\delta$. 
Then by Lemma \ref{firstlm},  $(v_{k_1})_{+} \leq \frac{1}{2}$ in $Q_{1/2}$ so that 
$$v\leq 1-\frac{1}{2^{k_1+1}}\leq 1-\frac{1}{2^{k_0+1}} \quad \mbox{ in } Q_{1/2},$$
So $\mu= \frac{1}{2^{k_0+1}}= \frac{1}{2^{|Q_2|\left(\frac{2C}{\delta|Q_1|}\right)^{4p+2}+2}}$ and by definition $\beta=\frac{1}{2^{k_0}}$ which implies that $v_k$ remains in $DG^{+}(\gamma_1,\gamma_2,\gamma_3,p)$ for $k\leq k_0$.
\end{proof}

 \subsection{Proof of Theorem \ref{t:Holder}.}
 In this subsection we provide a proof for Theorem~\ref{t:Holder}.
\begin{proof}[Proof of Theorem \ref{t:Holder})]
Let $u:Q_2\rightarrow\R$ be a function in $\cDG^+(\gamma_1, \gamma_2,\gamma_3, p)\cap \cDG^-(\gamma_1, \gamma_2,\gamma_3, p)$ then define $v:=\beta u$ that is in $\cDG^+(\gamma_1, \gamma_2,\beta\gamma_3, p)\cap \cDG^-(\gamma_1, \gamma_2,\beta\gamma_3, p)$.
Let $(t_0,x_0)\in Q_1$, $r\in (0,1/2)$ and for $n\in \mathbb{N}\setminus\{0\}$, define $$v_n(\tau,y)= \frac{2\theta^{1-n}}{\max(2,\underset{Q_{3/2}}{\mathrm{osc}} v )} v\left(t_0+\frac{\tau}{4^n},x_0+\frac{y}{2^n}\right)$$ and $$V_n= \frac{2}{\underset{Q_{1}}{\mathrm{osc}} \mbox{ } v_n } \left(v_n-\frac{\sup v_n +\inf v_n}{2}\right)$$ both in $  \cDG^+(\gamma_1, \gamma_2,\beta\gamma_3, p)\cap \cDG^-(\gamma_1, \gamma_2,\beta\gamma_3, p)$. Notice for $n\geq 2$ that $-1\leq V_{n-1} \leq 1$ and apply Lemma \ref{DGpLm5} to whichever of $V_{n-1}$ or $-V_{n-1}$ satisfies \eqref{hypvplm5} to get with $\theta:=1-\frac{\mu}{2}$
$$ \underset{Q_{1/2}}{\mathrm{osc }} v_{n-1} \leq \theta \max\left(2,\underset{Q_{1}}{\mathrm{osc }} \mbox{ } v_{n-1}\right)$$
so that by induction since $\underset{Q_{1/2}}{\mathrm{osc }} v_1 \leq 2$,  
$$\underset{Q_{1}}{\mathrm{osc }} \mbox{ } v_n=\underset{Q_{1/2}}{\mathrm{osc }} v_{n-1} \leq \theta \max(2,\underset{Q_{1}}{\mathrm{osc }} \mbox{ } v_{n-1})\leq 2\theta.$$
Choosing $\alpha \in (0,1)$ such that $\theta=\frac{1}{2^{\alpha}}$ and $n\in \mathbb{N}\setminus\{0\}$ such that $\frac{1}{2^{n+1}}\leq r<\frac{1}{2^n}$, deduce
\begin{align}
\label{ineq:reformulation}
\underset{Q_{r}(t_0,x_0)}{\mathrm{osc }} v &\leq \underset{Q_{\frac{1}{2^n}}(t_0,x_0)}{\mathrm{osc }}=\frac{\max(2,\underset{Q_{3/2}}{\mathrm{osc }} v)\theta^{n-1}}{2}\underset{Q_{1}}{\mathrm{osc }}\mbox{ } v_n 
 \leq \max(2,\underset{Q_{3/2}}{\mathrm{osc }} v) \theta^n 
\leq Cr^\alpha \left(\|v\|_{L^{2}(Q_{2})}+1\right),
 \end{align}
where we used Lemma \ref{firstlm} in the last inequality.
Conclude the H\"older estimate for $v$ between $z,z'\in Q_1$ using \eqref{ineq:reformulation} to intermediate points between $z$ and $z'$ at distance $r<\frac{1}{2}$ and a triangle inequality and deduce the H\"older estimate for $u$.
\end{proof} 
\appendix
\section{Proof of Corollary \ref{cor:LVIsol} and \ref{cor:Holder}} 
\label{A:sol in DG} 
 \begin{prop}
\label{sol in sub class}
Let $u$ be a sub-solution (resp. supersolution) of \eqref{eqpara} satisfying \eqref{hypcoeff}. Then there exist $\gamma_1,\gamma_2$ and $\gamma_3$ positive such that $u\in DG^{+}(\gamma_1,\gamma_2,\gamma_3,p)$ (resp. $u\in DG^{-}(\gamma_1,\gamma_2,\gamma_3,p)$).
Moreover if $u$ is a solution then there exist $\gamma_1, \gamma_2, \gamma_3$ and $p$ such that  $u\in DG^{+}(\gamma_1,\gamma_2,\gamma_3,p)\cap DG^{-}(\gamma_1,\gamma_2,\gamma_3,p)$. 
\end{prop}
 
\begin{proof}
Let us deal with the sub-solution case. The case of super-solution is very similar and the case of the solution is a combination a both cases. 
It is exactly deriving the energy estimates for the sub-solution $(u-k)_+$ of \eqref{eqpara} with the source term $g\mathds{1}_{(u-k)_+}$.
%We multiply the equation by $(u-k)+$ which is a subsolution of \eqref{eqpara} with the source term $s\mathds{1}_{(u-k)_+}$ and integrate in $(s,t)\times B_R$. 
Let us define $\varphi \in C^{\infty}_{c}(B_R)$ such that $0\leq \varphi \leq 1$, $|\nabla_x \varphi|\leq \frac{2}{R-r}$ and 
$$\varphi= \left\{\begin{array}{ll}
1 & \mbox{ in } B_r\\
0& \mbox{ ouside } B_R,
\end{array} \right.$$ 
and the sequence of functions 
$$\psi_{\varepsilon}(\tau)=\left\{\begin{array}{ll}
\frac{1}{\varepsilon}(\tau-s) & \mbox{ if }  \tau \in (s,s+\varepsilon)\\
1 & \mbox{ if } \tau \in (s+\varepsilon,t-\varepsilon)\\
-\frac{1}{\varepsilon}(\tau-t) & \mbox{ if }  \tau \in (t-\varepsilon,t).
\end{array} \right.$$
The idea is to use the test function $(\tau,x)\rightarrow (u-k)_+ (\tau,x)\varphi^2(x)$ which is not with compact support in time by first using the function $(\tau,x)\rightarrow (u-k)_+ (\tau,x)\psi_{\varepsilon}(\tau)\varphi^2(x)$ as a test function which is allowed by density arguments and then take $\varepsilon\rightarrow 0$. 
\begin{align*}
\displaystyle -\int_{Q} (u-k)_+\partial_{t}\left[(u-k)_+\psi_{\varepsilon}\varphi^2\right] & =\int_{Q} \partial_{t}(u-k)_+ (u-k)_+ \psi_{\varepsilon}\varphi^2 \\ & =\frac{1}{2}\int_{Q} \partial_{t}(u-k)_+^{2} \psi_{\varepsilon}\varphi^2 \\ &=- \frac{1}{2}\int_{Q} (u-k)_+^{2} \partial_t\psi_{\varepsilon}\varphi^2\\
&= -\frac{1}{2\epsilon} \displaystyle \int_{s}^{s+\varepsilon} \int_{\Omega}(u-k)_+^{2} \varphi + \frac{1}{2\epsilon} \displaystyle \int_{t-\varepsilon}^{t} \int_{\Omega} (u-k)_+^{2} \varphi^2
\end{align*}
The Lebesgue differentiation theorem for $\varepsilon \rightarrow 0$ gives
\begin{align}
\label{limit sous sol in DG}
\displaystyle -\int_{Q} (u-k)_+\partial_{t}\left[(u-k)_+\psi_{\varepsilon}\varphi^2\right] \rightarrow -\frac{1}{2} \int_{\Omega}(u-k)_+^{2}(s,.) \varphi^2 + \frac{1}{2} \displaystyle \int_{\Omega} (u-k)_+^{2}(t,.) \varphi^2.
\end{align}

For the other terms, since there is no derivative in time, $\psi_{\varepsilon}$ will be a common factor for each term. By using the dominated convergence theorem, when $\varepsilon \rightarrow 0$, $\psi_{\varepsilon} \rightarrow 1$ almost everywhere so at the limit, the other terms will be 
\begin{align*}
I:=\int_{Q} A\nabla_{x}(u-k)_+\cdot \nabla_{x}\left[(u-k)_+ \varphi^2 \mathds{1}_{(s,t)}\right]-\int_{Q} B\cdot \nabla_{x} (u-k)_+ (u-k)_+\varphi^2\mathds{1}_{(s,t)}\\ -\int_{Q} g(u-k)_+ \varphi^2\mathds{1}_{(s,t)}.
\end{align*}
Using a Young inequality, the inequality $\| g\|_{L^q} \leq 1$ and defining $p=\frac{q}{q-1}$,
\begin{align}
\label{ineg sous sol in DG}
I\geq & \lambda\int_{s}^t\int_{B_R} |\nabla_{x}(u-k)_+|^{2}\varphi^2 - 2\Lambda\int_{s}^t\int_{B_R} |\nabla_{x}(u-k)_+|\varphi (u-k)_+|\nabla_{x}\varphi|\nonumber \\ 
 & -\Lambda\int_{s}^t\int_{B_R} |\nabla_{x} (u-k)_+| (u-k)_+\varphi^2-\int_{s}^t\int_{B_R} |g|(u-k)_+ \varphi^2 \nonumber\\
\geq & \lambda\int_{s}^t\int_{B_R} |\nabla_{x}(u-k)_+|^{2}\varphi^2 -\frac{\lambda}{4}\int_{s}^t\int_{B_R} |\nabla_{x}(u-k)_+|^{2}\varphi^2 - \frac{4\Lambda^2}{\lambda}\int_{s}^t\int_{B_R} (u-k)_+^2|\nabla_{x}\varphi|^2 \nonumber\\  
& -\frac{\lambda}{4}\int_{s}^t\int_{B_R} |\nabla_{x}(u-k)_+|^{2}\varphi^2 -\frac{\Lambda^2}{\lambda}\int_{s}^t\int_{B_R} (u-k)_+^2|\nabla_{x}\varphi|^2 -\|g\|_{L^{q}}\left(\int_{s}^t\int_{B_R}(u-k)_+^{\frac{q}{q-1}} \varphi^2\right)^{\frac{q-1}{q}}\nonumber\\
\geq & \frac{\lambda}{2}\int_{s}^t\int_{B_r} |\nabla_{x}(u-k)_+|^{2} - \frac{5\Lambda^2}{\lambda(R-r)^2}\int_{s}^t\int_{B_R} (u-k)_+^2 -\|g\|_{L^{q}}\left(\int_{s}^t \int_{B_R}(u-k)_+^{\frac{q}{q-1}}\right)^{\frac{q-1}{q}}
\end{align}

Combining \eqref{limit sous sol in DG} and \eqref{ineg sous sol in DG} gives $u \in DG^{+}(\frac{\lambda}{2},\frac{5\Lambda^2}{\lambda}, \|g\|_{L^{q}}, \frac{q}{q-1})$. 
\end{proof}

\section{Tool for the $L^2-L^{\infty}$ estimate} 
\label{B:tool}

\begin{lemma}
\label{seq tends 0}
Let $(V_k)_{k\leq 0}$ be a sequence of real numbers such that for all $k\geq 1$,
\begin{equation}
\label{recurrence formula}
V_k\leq C^k V_{k-1}^{\alpha}
\end{equation}
where $\alpha>1$. 
Then for $V_0 < C^{-\frac{\alpha^2}{(\alpha-1)^2}}$, the sequence $(V_k)$ converges to $0$ when $k\rightarrow \infty$. 
\end{lemma}

\begin{proof}
By induction, deduce 
$$V_k\leq C^{k+(k-1)\alpha+\dots + 2\alpha^{k-2}+\alpha^{k-1}}V_0^{\alpha^k}=C^{S_k}V_0^{\alpha^k},$$
where $S_k= \sum_{i=0}^{k} i \alpha^{k-i}.$
Let us prove that for all $k\geq 1$,
\begin{equation}
\label{bound Sk}
S_k\leq \frac{\alpha^2}{(\alpha-1)^2}\alpha^{k-1}. 
\end{equation}
In fact, 
$$S_k= \alpha^{k-1} \sum_{i=0}^{k} i \left(\frac{1}{\alpha}\right)^{i-1}.$$
differentiate 
$$ \sum_{i=0}^{k}  X^{i} = \frac{1-X^{k+1}}{1-X},$$
to get
$$ \sum_{i=0}^{k}  i X^{i-1} = \frac{X^{k}(kX-(k+1))+1}{(1-X)^2}. $$
Since $\alpha>1$, deduce 
$$S_k\leq \alpha^{k-1} \frac{1}{\left(1-\frac{1}{\alpha}\right)^2},$$
which gives \eqref{bound Sk} and
$$ V_k\leq C^{\frac{\alpha^2}{(\alpha-1)^2}\alpha^{k-1}}V_0^{\alpha^k} \leq \left(C^{\frac{\alpha^2}{(\alpha-1)^2}}V_0\right)^{\alpha^k}.$$
which gives the limit for $V_0 < C^{-\frac{\alpha^2}{(\alpha-1)^2}}$.  
\end{proof}
 
 \vskip14mm
\noindent {\bf Acknowledgements.} 
The author acknowledges Cyril Imbert and Cl\'ement Mouhot for their expert advice throughout working on this paper. The author acknowledges the support of partial funding by the ERC grant MAFRAN 2017-2022.


\begin{thebibliography}{}
\bibitem{AD04}
Gabriel Acosta and Ricardo~G Dur{\'a}n.
\newblock An optimal {P}oincar{\'e} inequality in ${L}^1$ for convex domains.
\newblock {\em Proceedings of the american mathematical society}, pages
  195--202, 2004.

\bibitem{adams03}
Robert~A Adams and John~JF Fournier.
\newblock {\em Sobolev spaces}, volume 140.
\newblock Elsevier, 2003.

\bibitem{CCV}
Luis Caffarelli, Chi~Hin Chan, and Alexis Vasseur.
\newblock Regularity theory for parabolic nonlinear integral operators.
\newblock {\em J. Amer. Math. Soc.}, 24(3):849--869, 2011.

\bibitem{degiorgi}
Ennio De~Giorgi.
\newblock Sulla differenziabilit\`a e l'analiticit\`a delle estremali degli
  integrali multipli regolari.
\newblock {\em Mem. Accad. Sci. Torino. Cl. Sci. Fis. Mat. Nat. (3)}, 3:25--43,
  1957.

\bibitem{Dibenedetto89}
Emmanuele DiBenedetto.
\newblock Harnack estimates in certain function classes.
\newblock {\em Atti Semin. Mat. Fis. Univ. Modena}, 37(1):173--182, 1989.

\bibitem{DiBen2010}
Emmanuele DiBenedetto.
\newblock {\em Linear Elliptic Equations with Measurable Coefficients}, pages
  297--345.
\newblock Birkh{\"a}user Boston, Boston, 2010.

\bibitem{Gianazza_Vespri06}
Ugo Gianazza and Vincenzo Vespri.
\newblock Parabolic de Giorgi classes of order p and the Harnack inequality.
\newblock {\em Calculus of Variations and Partial Differential Equations},
  26(3):379--399, 2006.

\bibitem{Giaquinta_Gusti82}
Mariano Giaquinta and Enrico Giusti.
\newblock On the regularity of the minima of variational integrals.
\newblock {\em Acta Math.}, 148:31--46, 1982.

\bibitem{GIMV}
Fran\c{c}ois Golse, Cyril Imbert, Cl{\'e}ment Mouhot, and Alexis~F. Vasseur.
\newblock {Harnack inequality for kinetic Fokker-Planck equations with rough
  coefficients and application to the Landau equation}.
\newblock {\em Annali della Scuola Normale Superiore di Pisa, Classe di
  Scienze}, Vol. XIX, issue 1:PP. 253--295, 2019.

\bibitem{gue3}
Jessica Guerand.
\newblock Quantitative parabolic regularity \`a la {D}e~{G}iorgi.
\newblock {\em S\'eminaire Laurent Schwartz --- EDP et applications},
  2018-2019.

\bibitem{Heran22}
Andreas Her\'an.
\newblock H\"older continuity of parabolic quasi-minimizers on metric measure
  spaces.
\newblock {\em Journal of Differential Equations}, 341:208--262, 2022.

\bibitem{HouNiu}
Lingling Hou and Pengcheng Niu.
\newblock A Nash type result for divergence parabolic equation related to
  H\"ormander's vector fields.
\newblock {\em Journal of Partial Differential Equations}, 33(4):341--376,
  2020.

\bibitem{LSU}
O.~A. Ladyzhenskaya, V.~A. Solonnikov, and N.~N. Ural'tseva.
\newblock {\em Linear and quasilinear equations of parabolic type}.
\newblock Translated from the Russian by S. Smith. Translations of Mathematical
  Monographs, Vol. 23. American Mathematical Society, Providence, R.I., 1968.

\bibitem{LU1}
Olga~A. Ladyzhenskaya and Nina~N. Ural'tseva.
\newblock {\em Linear and quasilinear elliptic equations}.
\newblock Translated from the Russian by Scripta Technica, Inc. Translation
  editor: Leon Ehrenpreis. Academic Press, New York-London, 1968.

\bibitem{Liao21}
Naian Liao.
\newblock Remarks on parabolic de Giorgi classes.
\newblock {\em Annali di Matematica Pura ed Applicata (1923 -)},
  200(6):2361--2384, 2021.

\bibitem{Masson_Silijander13}
Mathias Masson and Juhana Siljander.
\newblock H{\"o}lder regularity for parabolic de Giorgi classes in metric
  measure spaces.
\newblock {\em Manuscripta Mathematica}, 142(1):187--214, 2013.

\bibitem{moser}
J\"urgen Moser.
\newblock A new proof of {D}e {G}iorgi's theorem concerning the regularity
  problem for elliptic differential equations.
\newblock {\em Comm. Pure Appl. Math.}, 13:457--468, 1960.

\bibitem{nash}
J.~Nash.
\newblock Continuity of solutions of parabolic and elliptic equations.
\newblock {\em Amer. J. Math.}, 80:931--954, 1958.

\bibitem{Vasseur2}
Alexis~F. Vasseur.
\newblock The {D}e {G}iorgi method for elliptic and parabolic equations and
  some applications.
\newblock In {\em Lectures on the analysis of nonlinear partial differential
  equations. {P}art 4}, volume~4 of {\em Morningside Lect. Math.}, pages
  195--222. Int. Press, Somerville, MA, 2016.

\bibitem{Wieser87}
Wilfried Wieser.
\newblock Parabolic q-minima and minimal solutions to variational flow.
\newblock {\em manuscripta mathematica}, 59(1):63--107, 1987.

          \end{thebibliography}
\end{document}